\documentclass[12pt]{article}
\usepackage{amssymb,amsthm,amsmath,latexsym}
\usepackage{url}
\newtheorem{thm}{Theorem}
\newtheorem{prop}[thm]{Proposition}

\theoremstyle{remark}
\newtheorem{rem}[thm]{Remark}
\newcommand{\FF}{\mathbb{F}}

\newcommand{\allone}{\mathbf{1}}

\DeclareMathOperator{\wt}{wt}
\DeclareMathOperator{\rank}{rank}

\begin{document}
\title{On extremal double circulant self-dual codes of lengths $90$--$96$}

\author{
T. Aaron Gulliver\thanks{Department of Electrical and Computer Engineering,
University of Victoria,
P.O. Box 1700, STN CSC, Victoria, BC,
Canada V8W 2Y2.
email: agullive@ece.uvic.ca}
\mbox{ and}
Masaaki Harada\thanks{
Research Center for Pure and Applied Mathematics,
Graduate School of Information Sciences,
Tohoku University, Sendai 980--8579, Japan.
email: mharada@m.tohoku.ac.jp.}
}
\date{}
\maketitle

\begin{abstract}
A classification of extremal double circulant self-dual codes
of lengths up to $88$ is known.
We give a classification of extremal double circulant self-dual
codes of lengths $90,92,94$ and $96$.
We also classify double circulant self-dual
codes with parameters $[90,45,14]$ and $[96,48,16]$.
In addition, we demonstrate that
no double circulant self-dual $[90,45,14]$ code
has an extremal self-dual neighbor, and
no double circulant self-dual $[96,45,16]$ code
has a self-dual neighbor with minimum
weight at least $18$.
\end{abstract}

%%%%%%%%%%%%%%%%%%%%%%%%%%%%%%
\section{Introduction}\label{sec:1}
A (binary) $[n,k]$ {\em code} $C$ is a $k$-dimensional vector subspace
of $\FF_2^n$,
where $\FF_2$ denotes the finite field of order $2$.
All codes in this note are binary.
The parameter $n$ is called the {\em length} of $C$.
The {\em weight} $\wt(x)$ of a vector $x \in \FF_2^n$ is
the number of non-zero components of $x$.
A vector of $C$ is a {\em codeword} of $C$.
The minimum non-zero weight of all codewords in $C$ is called
the {\em minimum weight} of $C$ and an $[n,k]$ code with minimum
weight $d$ is called an $[n,k,d]$ code.
The {\em weight enumerator} $W(C)$ of $C$ is given by
$W(C)= \sum_{i=0}^{n} A_i y^i$ where $A_i$ is the number of codewords of
weight $i$ in $C$.
Two codes are {\em equivalent} if one can be obtained from the other by a
permutation of coordinates.
The {\em dual} code $C^{\perp}$ of a code
$C$ of length $n$ is defined as
$
C^{\perp}=
\{x \in \FF_2^n \mid x \cdot y = 0 \text{ for all } y \in C\},
$
where $x \cdot y$ is the standard inner product.
A code $C$ is called {\em self-dual}
if $C=C^\perp$.
A self-dual code $C$ is called
{\em doubly even} and {\em singly even}
if all codewords have weight
$\equiv 0 \pmod 4$ and if some codeword has weight
$\equiv 2 \pmod 4$, respectively.

It was shown in~\cite{Mallows-Sloane} that
the minimum weight $d$ of a doubly even self-dual code
of length $n$ is bounded by $d \le 4[n/24]+4$.
We call a doubly even self-dual code meeting this upper
bound {\em extremal}.
% Currently, it is not known if an extremal
% doubly even self-dual $[96,48,20]$ code exists.
%
The largest possible minimum weights of
(singly even) self-dual codes of lengths up to $72$ are
given in~\cite[Table~I]{C-S}.
This work was extended to lengths up to $100$ in~\cite[Table~VI]{DGH}
(see~\cite[Table~2]{GO} and \cite[Table~I]{HKWY}).
According to~\cite{GH-DCC88},
in this note,
we say that a singly even self-dual code with
the largest possible minimum weight given in~\cite[Table~I]{C-S}
and~\cite[Table~VI]{DGH}
is {\em extremal}.
%%%%%%
The largest possible minimum weight
among singly even self-dual codes of lengths $90,92,94$ and $96$ is $16,16,18$
and $18$, respectively.
%%%%%%
Currently, it is not known if an extremal
self-dual $[90,45,16]$ code exists.
There is a self-dual $[90,45,14]$ code~\cite{DGH}.
Many extremal self-dual $[92,46,16]$ codes are known
(see~\cite{DRZ}, \cite{GPSA}, \cite{HN}, \cite{Y92} and
\cite{YW}, and references [5] and [10] in~\cite{HN}).
Currently, it is not known if an extremal
self-dual $[94,47,18]$ code exists.
There is a self-dual $[94,47,16]$ code~\cite{HKWY}.
Currently, it is not known if an extremal
doubly even self-dual $[96,48,20]$ code, or
an extremal singly even self-dual $[96,48,18]$ code, exists.
There is a doubly even self-dual $[96,48,16]$ code (see~\cite{DGH}), and
a singly even self-dual $[96,48,16]$ code~\cite{GO}.

Let $D_p$ and $D_b$ be codes with generator matrices of the form
\begin{equation}\label{pDCC}
\left(\begin{array}{ccccc}
{} & I_n  & {} & R & {} \\
\end{array}\right)
\end{equation}
and
\begin{equation}\label{bDCC}
\left(\begin{array}{ccccccccc}
{} & {} & {}      & {} & {} & 0     & 1  & \cdots &1  \\
{} & {} & {}      & {} & {} & 1     & {} & {}     &{} \\
{} & {} & I_{n+1} & {} & {} &\vdots & {} & R'     &{} \\
{} & {} & {}      & {} & {} & 1     & {} &{}      &{} \\
\end{array}\right),
\end{equation}
%% respectively, where $I$ is the identity matrix of order $n$
respectively, where $I_n$ is the identity matrix of order $n$,
and $R$ and $R'$ are $n \times n$ circulant matrices.
%%%%%%
An $n \times n$ circulant matrix has the form
\[
\left(
\begin{array}{ccccc}
r_0&r_1&r_2& \cdots &r_{n-1} \\
r_{n-1}&r_0&r_1& \cdots &r_{n-2} \\
\vdots &\vdots & \vdots && \vdots\\
r_1&r_2&r_3& \cdots&r_0
\end{array}
\right)
\]
so that each successive row is a cyclic shift of the previous one.
%%%%%%
The codes $D_p$ and $D_b$ are called {\em pure double circulant}
and {\em bordered double circulant}, respectively.
The two families are called double circulant codes.
Many of the best known self-dual codes are double circulant codes
(see~\cite{DGH}, \cite{GH-WE}, \cite{GH-DCC},
\cite{GH-DCC88} and \cite{HGK-DCC}).
Further, constructions exist that provide double circulant self-dual codes with the largest known
% minimum weight (see \cite{GS}, \cite{Moore}, \cite{TTHAA}).
minimum weight (see \cite{GS} and \cite{Moore}).
The bordered double circulant construction provides self-dual codes only when
%% the length is congruent to $\equiv 0 \pmod4$ for self-dual
the length is $\equiv 0 \pmod4$.
In addition, it is known~\cite{GH-WE} that
there is no bordered double circulant singly even
self-dual code of length $n \equiv 0 \pmod8$.

A classification of extremal double circulant self-dual codes
of lengths up to $88$ was given
in~\cite{GH-DCC}, \cite{GH-DCC88} and~\cite{HGK-DCC}.
In this note, this work is extended to length $96$.
Our exhaustive search shows that there is no extremal double
circulant self-dual $[90,45,16]$ code.
We also give a classification of double circulant self-dual
$[90,45,14]$ codes.
In addition, we demonstrate that
every double circulant self-dual $[90,45,14]$ code
has no extremal self-dual neighbor.
We give a classification of extremal double circulant self-dual
codes of length $92$.
Our exhaustive search shows that there is no double circulant self-dual
$[94,47,d]$ code with $d \ge 16$.
We give a classification of double circulant self-dual
$[96,48,d]$ codes with $d \ge 16$.
In addition, we demonstrate that
every double circulant self-dual $[96,48,16]$ code
has no self-dual neighbor with minimum weight at least $18$.

%%%%%%%%%%%%%%%%%%%%
\section{Double circulant self-dual $[90,45,d]$ codes with $d \in \{14,16\}$}
\label{sec:90}

Using an approach similar to that given in~\cite{GH-DCC},
\cite{GH-DCC88} and~\cite{HGK-DCC},
our exhaustive search
found all distinct double circulant
self-dual $[90,45,d]$ codes with $d \ge 14$.
This was done by considering all
$45 \times 45$ orthogonal circulant matrices satisfying
the condition that the weight of the first row is congruent to $1 \pmod 4$
and the weight is greater than or equal to $d-1$.
Since a cyclic shift of the first row of some
codes defines an equivalent code,
the elimination of cyclic shifts
substantially reduces the number of codes
which must be checked further for equivalence to
complete the classification.
It is useful to use the fact that
self-dual codes with generator matrices of the form
$
\left(\begin{array}{ccccc}
I_{45}  & R  \\
\end{array}\right)
$
and
$
\left(\begin{array}{ccccc}
I_{45}  & R^T  \\
\end{array}\right)
$ are equivalent,
where $R^T$ denotes the transpose of $R$.
{\sc Magma}~\cite{Magma} was employed to determine code
equivalence and complete the classification.
Then we have the following results.

\begin{prop}
There is no extremal double
circulant self-dual code of length $90$.
\end{prop}

\begin{prop}
There are $716$ inequivalent double circulant self-dual
$[90,45,14]$ codes.
\end{prop}

The first rows of $R$ in the generator matrices
$
\left(\begin{array}{ccccc}
I_{45}  & R  \\
\end{array}\right)
$
of the $716$ codes can be obtained from
\url{http://www.math.is.tohoku.ac.jp/~mharada/Paper/DCC90.txt}.
We verified by {\sc Magma}~\cite{Magma} that each of the
$716$ codes has an automorphism group of order $90$.

% We determine the possible weight enumerators
% of self-dual $[90,45,14]$ codes and their shadows
% (see \cite{C-S} for the definition of shadows).
% For a detailed description of how this is accomplished,
% see~\cite[Theorem~5]{C-S}.
% The possible weight enumerators $W_i$ (resp.\ $S_i$) $(i=1,2,3)$
% of self-dual $[90,45,14]$ codes (resp.\ their shadows)
% are as follows:
We determined the possible weight enumerators
of self-dual $[90,45,14]$ codes.
For a detailed description of how this is accomplished,
see~\cite[Theorem~5]{C-S}.
The possible weight enumerators of self-dual $[90,45,14]$ codes
and the shadows are as follows
\begin{align*}
&
1
+(14040 + a)y^{14}
+(51300 + 3 a + 8 b)y^{16}
+(69920 -  11 a - 24 b + 512 c)y^{18}
\\&
+(2355624 - 41 a - 80 b -  4608 c + 32768 d)y^{20}
\\&
+(30913560 + 49 a + 304 b +  13824 c - 491520 d - 2097152 e)y^{22}
+ \cdots,
\\&
e y
+(d - 22 e)y^5
+(- c - 20 d + 231 e)y^9
+(b + 18 c + 190 d - 1540 e)y^{13}
\\&
+(- 8 a - 16 b - 153 c -  1140 d + 7315 e)y^{17}
% \\&
% +(20820480 + 112 a + 120 b + 816 c + 4845 d - 26334 e)y^{21}
+ \cdots,
\end{align*}
respectively, where $a,b,c,d,e$ are integers.
It is easy to see that
the number of codewords of weights $14,16$ in the code
and the number of vectors of weights $1,5,9$ in the shadow
uniquely determine the weight enumerator.
By calculating these numbers,
we verified that
the $716$ codes have $100$ distinct weight enumerators.
This was done using {\sc Magma}~\cite{Magma}.
The $100$ weight enumerators have $c=d=e=0$,
where $(a,b)$ are listed in Table~\ref{Tab:WD}.
For each pair $(a,b)$,
the number $N_{(a,b)}$ of codes with the weight enumerator
is also listed in Table~\ref{Tab:WD}.

%%%%%%%%%%%%%%%%%%%%%%%%%%%%%%%
\begin{table}[thbp]
\caption{Weight enumerators of double circulant self-dual $[90,45,14]$ codes}
\label{Tab:WD}
\begin{center}
%{\small
{\footnotesize
%{\scriptsize
\begin{tabular}{c|c|c|c|c|c}
\noalign{\hrule height0.8pt}
$(a,b)$ &  $N_{(a,b)}$ &
$(a,b)$ &  $N_{(a,b)}$ &
$(a,b)$ &  $N_{(a,b)}$ \\
\hline
$(-12555, 0)$  & 1&
$(-12555, 90)$ & 1&
$(-12600, 180)$& 1\\
$(-12735, 0)$  & 1&
$(-12780, 0)$  & 1&
$(-12825, 0)$  & 3\\
$(-12825, 90)$ & 2&
$(-12870, 0)$  & 5&
$(-12870, 180)$& 1\\
$(-12915, 0)$  & 3&
$(-12915, 180)$& 1&
$(-12915, 90)$ & 4\\
$(-12960, 0)$  & 6&
$(-12960, 90)$ & 9&
$(-13005, 0)$  & 6\\
$(-13005, 180)$& 1&
$(-13005, 90)$ &10&
$(-13050, 0)$  &10\\
$(-13050, 180)$& 1&
$(-13050, 270)$& 1&
$(-13050, 90)$ & 8\\
$(-13095, 0)$  & 8&
$(-13095, 180)$& 4&
$(-13095, 270)$& 1\\
$(-13095, 90)$ & 6&
$(-13140, 0)$  &17&
$(-13140, 180)$& 7\\
$(-13140, 270)$& 1&
$(-13140, 90)$ &13&
$(-13185, 0)$  &10\\
$(-13185, 180)$& 5&
$(-13185, 270)$& 1&
$(-13185, 360)$& 1\\
$(-13185, 90)$ &16&
$(-13230, 0)$  & 6&
$(-13230, 180)$&10\\
$(-13230, 270)$& 4&
$(-13230, 90)$ &27&
$(-13275, 0)$  &21\\
$(-13275, 180)$&15&
$(-13275, 270)$& 1&
$(-13275, 360)$& 3\\
$(-13275, 90)$ &16&
$(-13320, 0)$  &11&
$(-13320, 180)$&20\\
$(-13320, 270)$& 5&
$(-13320, 450)$& 1&
$(-13320, 90)$ &19\\
$(-13365, 0)$  &13&
$(-13365, 180)$&19&
$(-13365, 270)$& 7\\
$(-13365, 360)$& 1&
$(-13365, 90)$ &13&
$(-13410, 0)$  & 8\\
$(-13410, 180)$&22&
$(-13410, 270)$& 6&
$(-13410, 360)$& 2\\
$(-13410, 90)$ &23&
$(-13455, 0)$  &11&
$(-13455, 180)$&18\\
$(-13455, 270)$&13&
$(-13455, 360)$& 1&
$(-13455, 90)$ &33\\
$(-13500, 0)$  & 8&
$(-13500, 180)$&17&
$(-13500, 270)$& 9\\
$(-13500, 360)$& 3&
$(-13500, 90)$ &23&
$(-13545, 0)$  & 4\\
$(-13545, 180)$&19&
$(-13545, 270)$& 8&
$(-13545, 450)$& 1\\
$(-13545, 90)$ &18&
$(-13590, 0)$  & 3&
$(-13590, 180)$& 8\\
$(-13590, 270)$& 9&
$(-13590, 360)$& 2&
$(-13590, 450)$& 2\\
$(-13590, 90)$ & 9&
$(-13635, 0)$  & 2&
$(-13635, 180)$& 7\\
$(-13635, 270)$& 5&
$(-13635, 360)$& 2&
$(-13635, 450)$& 2\\
$(-13635, 90)$ & 6&
$(-13680, 180)$& 9&
$(-13680, 270)$& 4\\
$(-13680, 360)$& 1&
$(-13680, 90)$ & 2&
$(-13725, 180)$& 2\\
$(-13725, 270)$& 2&
$(-13725, 360)$& 2&
$(-13725, 90)$ & 2\\
$(-13770, 180)$& 1&
$(-13770, 270)$& 2&
$(-13815, 180)$& 1\\
$(-13815, 270)$& 2&
$(-13815, 360)$& 1&
$(-13815, 90)$ & 1\\
$(-13905, 360)$& 2& & & \\
\noalign{\hrule height0.8pt}
\end{tabular}
}
\end{center}
\end{table}
%%%%%%%%%%%%%%%%%%%%%%%%%%%%%%%

%%%%%%%%%%%%%%%%%%%%%
\section{Neighbors of double circulant self-dual $[90,45,14]$ codes}
\label{sec:90nei}

Two self-dual codes $C$ and $C'$ of length $n$
are said to be {\em neighbors} if $\dim(C \cap C')=n/2-1$.
We give some observations from~\cite{CHK64} on self-dual codes constructed by neighbors.
Let $C$ be a self-dual $[n,n/2,d]$ code.
Let $M$ be a matrix whose rows are the codewords of weight $d$ in $C$.
Suppose that there is a vector $x\in \mathbb{F}_2^n$ such that
\begin{equation}\label{Eq:N}
M x^T = \allone^T,
\end{equation}
where
$\allone$ is the all-one vector.
Set $C_0=\langle x\rangle^\perp\cap C$, where $\langle x\rangle$ denotes
the code generated by $x$.
Then $C_0$ is a subcode of index $2$ in $C$.
If the weight of $x$ is even, then
we have two self-dual neighbors $\langle C_0,x\rangle$ and $\langle
C_0,x+y\rangle$ of $C$
for some $y\in C\setminus C_0$,
which do not contain any codewords of weight $d$ in $C$,
where $\langle C,x\rangle = C \cup (x+C)$.
When $C$ has a self-dual $[n,n/2,d']$ neighbor $C'$ with $d'\ge d+2$,
\eqref{Eq:N} has a solution $x$ and we can obtain $C'$ in this way.
If $\rank M < \rank ( M\ \allone^T)$,
then $C$ has no self-dual $[n,n/2,d']$ neighbor $C'$ with $d'\ge d+2$.
%% Moreover, if $C$ is doubly even then one of the two self-dual
%% neighbors is doubly even and the other is singly even.
If $\rank M=t$, then we have at most
$2\times 2^{n/2-t}$ self-dual neighbors of $C$.
Furthermore, if the subcode
generated by the codewords of weight $d$ in $C$
contains $\allone$,
then $C$ has exactly $2\times 2^{n/2-t}$ self-dual neighbors.
When $C$ has a self-dual $[n,n/2,d']$ neighbor $C'$ with $d'\ge d+2$,
\eqref{Eq:N} has a solution $x$ and we can obtain $C'$ in this way.

We verified by {\sc Magma}~\cite{Magma} that
\[
(\rank M,\rank ( M \ \allone^T))=(43,43),
\]
for one of the $716$ double circulant self-dual $[90,45,14]$ codes and
\[
(\rank M,\rank ( M \ \allone^T))=(45,45),
\]
for the remaining $715$ codes.
% In addition, using a method similar to that given in~\cite{CHK64},
In addition, using the above method,
we verified by {\sc Magma}~\cite{Magma} that
the self-dual neighbors constructed by the above argument
have minimum weight at most $14$.
Hence, we have the following result.

\begin{prop}
No double circulant self-dual $[90,45,14]$ code
has an extremal self-dual neighbor of length $90$.
\end{prop}

It is still an open problem whether
an extremal self-dual code of length $90$ exists.

%%%%%%%%%%%%%%%%%%%%%
\section{Extremal double circulant self-dual codes of length $92$}

Using a method similar to that given in Section~\ref{sec:90},
our exhaustive search found all distinct extremal pure and bordered double circulant
self-dual codes of length $92$.
Then we have the following results.

\begin{prop}
There is no extremal pure double circulant
self-dual code of length $92$.
\end{prop}

\begin{rem}
Alfred Wassermann in a private communication indicated
that there is no extremal pure double circulant
self-dual code of length $92$, which provides an independent confirmation of our results.
\end{rem}

\begin{prop}
There are $158$ inequivalent extremal bordered double circulant
self-dual codes of length $92$.
\end{prop}

We denote the
$158$ inequivalent extremal bordered double circulant
self-dual codes of length $92$ by $B_{92,i}$ $(i=1,2,\ldots,158)$.
For the codes $B_{92,i}$ $(i=1,2,\ldots,158)$,
the first rows $r$ of $R'$ in~\eqref{bDCC} are listed in
Table~\ref{Tab:row92}.
In the table, the rows are written in octal using $0=(000)$,
$1=(001),\ldots,6=(110)$ and $7=(111)$.
We verified by {\sc Magma}~\cite{Magma} that $B_{92,i}$
has an automorphism group of order $90$ for $i=1,2,\ldots,158$.

%%%%%%%%%%%%%%%%%%%%%%%%%%%%%%%
\begin{table}[thb]
\caption{Weight enumerators of $B_{92,i}$ $(i=1,2,\ldots,158)$}
\label{Tab:WD92}
\begin{center}
%{\small
{\footnotesize
%{\scriptsize
\begin{tabular}{c|l}
\noalign{\hrule height0.8pt}
$\beta$ &  \multicolumn{1}{c}{$i$} \\
\hline
1527&
81, 140
\\
1572&
57, 64, 102
\\
1617&
52, 55, 77, 106
\\
1662&
6, 14, 18, 33, 56, 87, 121, 151
\\
1707&
3, 36, 38, 50, 69, 99, 101, 109, 111, 123, 143
\\
1752&
5, 21, 40, 49, 63, 116, 125, 127
\\
1797&
15, 27, 32, 46, 89, 95, 105, 138, 141, 147, 152, 153, 156
\\
1842&
1, 8, 10, 17, 22, 66, 72, 85, 90, 97, 108
\\
1887&
13, 26, 39, 41, 44, 48, 58, 62, 74, 84, 91, 103, 110, 112, 113, 119, 130,
\\&
136, 139, 154, 155
\\
1932&
16, 51, 80, 98, 131, 134
\\
1977&
2, 23, 24, 47, 53, 59, 61, 86, 120, 126
\\
2022&
28, 31, 37, 60, 67, 79, 82, 88, 92, 114, 117, 118, 128, 137, 144, 146, 150
\\
2067&
4, 7, 19, 35, 100, 107, 157, 158
\\
2112&
20, 65, 76, 94, 96, 104, 129, 132, 142, 145
\\
2157&
11, 34, 45, 68, 70, 133
\\
2202&
25, 54, 71, 75, 83, 149
\\
2247&
30, 42, 115
\\
2292&
12, 29, 43, 93, 122
\\
2337&
73, 124, 135
\\
2382&
78
\\
2427&
148
\\
2607&
9
\\
\noalign{\hrule height0.8pt}
\end{tabular}
}
\end{center}
\end{table}
%%%%%%%%%%%%%%%%%%%%%%%%%%%%%%%

%%%%%%%%%%%%%%%%%%%%%%%%%%%%%%%
\begin{table}[p]
\caption{First rows of $R'$ in~\eqref{bDCC} for $B_{92,i}$
$(i=1,2,\ldots,158)$}
\label{Tab:row92}
\begin{center}
%{\small
%{\footnotesize
{\scriptsize
\begin{tabular}{c|c|c|c|c|c}
\noalign{\hrule height0.8pt}
$i$ &  $r$ & $i$ &  $r$ & $i$ &  $r$ \\
\hline
  1 & 045722771307000 &
  2 & 046354263735000 &
  3 & 054130272607000 \\
  4 & 102155447541000 &
  5 & 104521145473000 &
  6 & 110545607071000 \\
  7 & 111222513531000 &
  8 & 115060555603000 &
  9 & 126644436177000 \\
 10 & 130023052373000 &
 11 & 130115321127000 &
 12 & 141232724213000 \\
 13 & 165172671757000 &
 14 & 206201322771000 &
 15 & 207022470467000 \\
 16 & 221254213777000 &
 17 & 233413676413000 &
 18 & 236265461527000 \\
 19 & 243271301677000 &
 20 & 246165167155000 &
 21 & 263645737137000 \\
 22 & 265543417117000 &
 23 & 271737137037000 &
 24 & 304151364577000 \\
 25 & 306307103017000 &
 26 & 330141216433000 &
 27 & 341136314537000 \\
 28 & 367621175177000 &
 29 & 406233754353000 &
 30 & 407777023131000 \\
 31 & 436453744513000 &
 32 & 442536745265000 &
 33 & 453136757327000 \\
 34 & 506762231273000 &
 35 & 522423127177000 &
 36 & 533555410563000 \\
 37 & 536456502533000 &
 38 & 545721744703000 &
 39 & 552732211353000 \\
 40 & 577604234513000 &
 41 & 605303705657000 &
 42 & 616171605527000 \\
 43 & 616352763673000 &
 44 & 645661771573000 &
 45 & 652610723547000 \\
 46 & 041045754474740 &
 47 & 041074234673640 &
 48 & 043135457542240 \\
 49 & 043250274476540 &
 50 & 043574372741740 &
 51 & 043603253362640 \\
 52 & 044104457046640 &
 53 & 044226336701740 &
 54 & 044671546702540 \\
 55 & 044736177354640 &
 56 & 046241222776640 &
 57 & 047663617660740 \\
 58 & 050576621023740 &
 59 & 050751623626240 &
 60 & 052354251553140 \\
 61 & 053047243053740 &
 62 & 053452447124740 &
 63 & 054716632374740 \\
 64 & 055727742734140 &
 65 & 057165135375140 &
 66 & 060736704331140 \\
 67 & 061065253646540 &
 68 & 063057644302740 &
 69 & 063731763152340 \\
 70 & 065235232367740 &
 71 & 066764121737540 &
 72 & 072237363775740 \\
 73 & 101742440560540 &
 74 & 103257370547740 &
 75 & 104571361141740 \\
 76 & 105644361251740 &
 77 & 107752466111140 &
 78 & 112175633752540 \\
 79 & 112453746103640 &
 80 & 113343466667340 &
 81 & 114740315712340 \\
 82 & 115256766234740 &
 83 & 115303767255340 &
 84 & 115331337561640 \\
 85 & 116226336753640 &
 86 & 116277617462540 &
 87 & 123320123173340 \\
 88 & 123663146657540 &
 89 & 123763453707140 &
 90 & 127654632533640 \\
 91 & 127737665370740 &
 92 & 131553671516340 &
 93 & 132721322740340 \\
 94 & 134537632035740 &
 95 & 136446677675740 &
 96 & 140352750117540 \\
 97 & 141147475510340 &
 98 & 142376737627740 &
 99 & 145756370077140 \\
100 & 146514734137740 &
101 & 146542652434540 &
102 & 147307747376740 \\
103 & 151056130545740 &
104 & 153556753442340 &
105 & 153743314476340 \\
106 & 154356132761740 &
107 & 155237453760340 &
108 & 162603763561740 \\
109 & 162755377664740 &
110 & 163131275323740 &
111 & 203315731776440 \\
112 & 204336305345340 &
113 & 205721743736640 &
114 & 206127347363740 \\
115 & 216236374745540 &
116 & 216712337262740 &
117 & 216776346322340 \\
118 & 223601644714740 &
119 & 225103675656740 &
120 & 225756665264340 \\
121 & 227266265646740 &
122 & 227656146713640 &
123 & 231235753751440 \\
124 & 233466766660640 &
125 & 233475224764740 &
126 & 235161677207340 \\
127 & 236612727317440 &
128 & 236657266701540 &
129 & 244353765463340 \\
130 & 246531765347240 &
131 & 257573767412740 &
132 & 261574375164340 \\
133 & 265576361776640 &
134 & 272764241653740 &
135 & 275373174710140 \\
136 & 277333517335540 &
137 & 277444625674540 &
138 & 277464773666340 \\
139 & 277757446333340 &
140 & 306156645077740 &
141 & 306533337306340 \\
142 & 310773753761740 &
143 & 313675510752340 &
144 & 317071170770740 \\
145 & 324761561777740 &
146 & 357234774374740 &
147 & 463637676050640 \\
148 & 512662622756740 &
149 & 516703753357740 &
150 & 535715465737340 \\
151 & 537761713627540 &
152 & 547776766153140 &
153 & 626053773755740 \\
154 & 656707543160740 &
155 & 043557136776166 &
156 & 047172572571772 \\
157 & 051150777762736 &
158 & 066775577156766 & &\\
\noalign{\hrule height0.8pt}
\end{tabular}
}
\end{center}
\end{table}
%%%%%%%%%%%%%%%%%%%%%%%%%%%%%%%

The possible weight enumerators of extremal self-dual codes
of length $92$ are given in~\cite{DGH} as follows
\begin{align*}
W_{92,1} =&  1 + (4 \beta + 4692)y^{16} + (174800 - 8 \beta + 256
\alpha)y^{18} \\
& + (-2048 \alpha + 2425488 - 52 \beta)y^{20} +
\cdots, \\
W_{92,2} =&  1 + (4 \beta + 4692)y^{16} +(174800 -8\beta +256 \alpha)y^{18}\\
& + (-2048 \alpha + 2441872 - 52 \beta)y^{20} +
\cdots, \\
W_{92,3} =&  1 + (4 \beta + 4692)y^{16} + (121296 - 8 \beta)y^{18}
\\& +(3213968 -52 \beta)y^{20} + \cdots,
\end{align*}
where $\alpha, \beta$ are integers.
By calculating the numbers of codewords of weights $16,18,20$
in the codes,
we verified that $B_{92,i}$ has weight enumerator
$W_{92,3}$, where $i$ and $\beta$ in $W_{92,3}$ are listed in
Table~\ref{Tab:WD92}.

%%%%%%%%%%%%%%%%%%%%
\section{Double circulant self-dual $[94,47,d]$ codes with $d \ge 16$}
\label{sec:94}

As mentioned in Section~\ref{sec:1},
it is currently not known if an extremal
self-dual code of length $94$ exists.
There is a self-dual $[94,47,16]$ code~\cite{HKWY}.

Using a method similar to that given in Section~\ref{sec:90},
our exhaustive search found no double circulant
self-dual $[94,46,d]$ code with $d \ge 16$.
Then we have the following result.

\begin{prop}
There is no double circulant
self-dual code of length $94$ and minimum weight $d \ge 16$.
\end{prop}

%It is still an open problem to determine whether there is
%an extremal self-dual code of length $94$.

%%%%%%%%%%%%%%%%%%%%
\section{Double circulant self-dual $[96,48,d]$ codes with $d \ge 16$}
\label{sec:96}

As described in Section~\ref{sec:1},
it is currently not known if an extremal
doubly even self-dual $[96,48,20]$ code exists, or
if an extremal singly even self-dual $[96,48,18]$ code exists.
There is a doubly even self-dual $[96,48,16]$ code~\cite{DGH}, and
a singly even self-dual $[96,48,16]$ code~\cite{GO}.

Using a method similar to that given in Section~\ref{sec:90},
our exhaustive search found all distinct pure double circulant
self-dual $[96,48,d]$ codes with $d \ge 16$ and
all distinct bordered double circulant doubly even
self-dual $[96,48,d]$ codes with $d \ge 16$.
Then we have the following result.

\begin{prop}
There is no extremal double
circulant doubly even self-dual code of length $96$.
There is no extremal double
circulant singly even self-dual code of length $96$.
\end{prop}
\begin{rem}
Alfred Wassermann in a private communication indicated
that there is no double circulant
self-dual code of length $96$ and minimum weight $d \ge 18$,
which provides an independent confirmation of our results.
\end{rem}

\begin{prop}
There are $49$ inequivalent pure double circulant
singly even self-dual $[96,48,16]$ codes.
There are $4565$ inequivalent pure double circulant
doubly even self-dual $[96,48,16]$ codes.
There are $1532$ inequivalent bordered double circulant
doubly even self-dual $[96,48,16]$ codes.
\end{prop}

%%%%%%%%%%%%%%%%%%%%%%%%%%%%%%%
\begin{table}[thb]
\caption{First rows of $R$ in~\eqref{pDCC} for $C_{96,i}$ $(i=1,2,\ldots,49)$}
\label{Tab:row96S}
\begin{center}
%{\small
%{\footnotesize
{\scriptsize
\begin{tabular}{c|c|c|c|c|c}
\noalign{\hrule height0.8pt}
$i$ &  $r$ & $(a,b,c,d)$ & $i$ &  $r$ & $(a,b,c,d)$ \\
\hline
 1 & 5532465545470000 & $(9798,0,0,0) $&
 2 & 1011116717627400 & $(10050,0,0,0)$\\
 3 & 2117213667133520 & $(10164,0,0,0)$&
 4 & 5160450553527400 & $(10416,0,0,0)$\\
 5 & 0411642402747400 & $(10422,0,0,0)$&
 6 & 1110737636054400 & $(10434,0,0,0)$\\
 7 & 2730315332407400 & $(10566,0,0,0)$&
 8 & 4127775466731720 & $(10566,0,0,0)$\\
 9 & 5247741422235400 & $(10740,0,0,0)$&
10 & 1104701417751460 & $(10854,0,0,0)$\\
11 & 1334257665167760 & $(10980,0,0,0)$&
12 & 1551523722207400 & $(11154,0,0,0)$\\
13 & 1072513135756620 & $(11364,0,0,0)$&
14 & 5302720447547400 & $(11508,0,0,0)$\\
15 & 1115027566566720 & $(11820,0,0,0)$&
16 & 1176414666173320 & $(12108,0,0,0)$\\
17 & 1252510325477400 & $(12180,0,0,0)$ &
18 & 0707334570645560 & $(9618,-48,0,0)$\\
19 & 0536450432504760 & $(10326,-48,0,0)$&
20 & 2727311567536720 & $(10326,-48,0,0)$\\
21 & 4260735067342400 & $(10422,-48,0,0)$&
22 & 5720417224633400 & $(10434,-48,0,0)$\\
23 & 0465637224357620 & $(10566,-48,0,0)$&
24 & 0447671345066400 & $(11124,-48,0,0)$\\
25 & 0644667174474660 & $(11844,-48,0,0)$&
26 & 1667375134475360 & $(11994,-48,0,0)$\\
27 & 1233543431133400 & $(10458,-96,0,0)$&
28 & 5772526161347400 & $(10806,-96,0,0)$\\
29 & 4072735065262400 & $(11190,-96,0,0)$&
30 & 1077057777245360 & $(11634,-96,0,0)$\\
31 & 1473646640067400 & $(11670,-96,0,0)$&
32 & 7242336777667400 & $(11748,-96,0,0)$\\
33 & 0411474700534400 & $(11796,-96,0,0)$&
34 & 7005761177137400 & $(11940,-96,0,0)$\\
35 & 0420777500236160 & $(11940,-96,0,0)$&
36 & 0407113175431520 & $(11952,-96,0,0)$\\
37 & 0603237114035560 & $(12390,-96,0,0)$&
38 & 2577350620527260 & $(12852,-96,0,0)$\\
39 & 0670641356075760 & $(10818,-144,0,0)$&
40 & 1462237456233660 & $(11616,-144,0,0)$\\
41 & 5176656144756400 & $(12090,-144,0,0)$&
42 & 0463117521602660 & $(12132,-144,0,0)$\\
43 & 0411766016336120 & $(12198,-144,0,0)$&
44 & 1237215132353660 & $(12384,-144,0,0)$\\
45 & 2271227740255400 & $(12690,-144,0,0)$&
46 & 2114462227575400 & $(13050,-144,0,0)$\\
47 & 1176617376233720 & $(13218,-144,0,0)$&
48 & 0477501403733400 & $(12282,-192,0,0)$\\
49 & 5764337750370000 & $(14124,-288,0,0)$& \\
\noalign{\hrule height0.8pt}
\end{tabular}
}
\end{center}
\end{table}
%%%%%%%%%%%%%%%%%%%%%%%%%%%%%%%

We denote the $49$ inequivalent pure double circulant
singly even self-dual $[96,48,16]$ codes by
$C_{96,i}$ $(i=1,2,\ldots,49)$.
For these codes,
the first rows $r$ of $R$ in~\eqref{pDCC} are listed in
Table~\ref{Tab:row96S}.
In the table, the rows are written in octal using $0=(000)$,
$1=(001),\ldots,6=(110)$ and $7=(111)$.

The possible weight enumerators of singly even self-dual
$[96,48,d]$ codes with $d \ge 16$ and their shadows
(see~\cite{C-S} for the definition) are
\begin{align*}
&
1
+ (- 5814+ a)y^{16}
+ (97280+ 64b)y^{18}
+ (1694208- 16a- 384b+4096c)y^{20}
\\ &
+ (18969600+ 192b- 49152c- 262144d)y^{22}
\\ &
+ (184315200+ 120a+ 3328b+ 237568c+ 4718592d)y^{24}
+ \cdots,\\
&
dy^4
+ (c- 22d)y^8
+ (- b- 20c+ 231d)y^{12}
%\\ &
+ (a+ 18b+ 190c- 1540d)y^{16}
\\ &
+ (3231744- 16a- 153b- 1140c+ 7315d)y^{20}
\\ &
+ (369664000+ 120a+ 816b+ 4845c- 26334d)y^{24}
+ \cdots,
\end{align*}
respectively, where $a,b,c,d$ are integers.
It is easy to see that
the number of codewords of weight $16$ in the code
and the numbers of vectors of weights $4,8,12$ in the shadow
uniquely determine the weight enumerator.
By calculating these numbers,
we determined the weight enumerators of the codes $C_{96,i}$.
This was done using {\sc Magma}~\cite{Magma}.
We display in Table~\ref{Tab:row96S}
$(a,b,c,d)$ for the weight enumerators of
the codes $C_{96,i}$.
We verified by {\sc Magma}~\cite{Magma} that each of
the $49$ codes has an automorphism group of order $96$.

We denote the $4565$ inequivalent pure double circulant
doubly even self-dual $[96,48,16]$ codes by
$P_{96,i}$ $(i=1,2,\ldots,4565)$.
The first rows $r$ of $R$ in~\eqref{pDCC} can be obtained from
\url{http://www.math.is.tohoku.ac.jp/~mharada/Paper/DCCp96.txt}.
By the Gleason theorem (see~\cite{Mallows-Sloane}),
the possible weight enumerators of doubly even self-dual
$[96,48,d]$ codes with $d \ge 16$ are
\begin{align*}
&
1
+ a y^{16}
+ (3217056 - 16 a )y^{20}
+ (369844880 + 120 a )y^{24}
\\ &
+ (18642839520 - 560 a )y^{28}
+ (422069980215 + 1820 a )y^{32}
+ \cdots,
\end{align*}
where $a$ is an integer with $0 \le a \le 201066$.
By calculating the number of codewords of weight $16$,
we verified that
the $4565$ codes have $614$ distinct weight enumerators.
The numbers $a$ in the weight enumerators
are listed in Table~\ref{Tab:WD96p}.
% The smallest and largest integers $a$ among the weight enumerators
% of $P_{96,i}$ $(i=1,2,\ldots,4565)$ are $5808$ and $14718$,
% respectively.
We verified by {\sc Magma}~\cite{Magma} that
$4530$, $34$ and $1$ of the $4565$ codes have
automorphism groups of orders $96$, $192$ and
$89280$, respectively.
For the unique code with an automorphism group of order $89280$,
the first row $r$ of $R$ in~\eqref{pDCC} is
\[
(010001101111001001011111101110100100101111110000).
\]

%%%%%%%%%%%%%%%%%%%%%%%%%%%%%%%
\begin{table}[thbp]
\caption{Weight enumerators of pure double circulant 
doubly even self-dual $[96,48,16]$ codes}
\label{Tab:WD96p}
\begin{center}
%{\small
%{\footnotesize
{\scriptsize
\begin{tabular}{l}
\noalign{\hrule height0.8pt}
\multicolumn{1}{c}{$a$} \\
\hline
5808, 6222, 6444, 6732, 6780, 6816, 6828, 6876, 6906, 6924, 6942, 6972, 6990, 
7002, 7008, 7068, 7092, 
\\
7098, 7104, 7116, 7152, 7182, 7236, 7278, 7290, 7296, 7308, 7332, 7338, 7374, 
7380, 7392, 7398, 7404, 
\\
7422, 7440, 7446, 7452, 7470, 7476, 7482, 7488, 7500, 7518, 7524, 7548, 7578, 
7584, 7596, 7614, 7620, 
\\
7626, 7632, 7638, 7644, 7662, 7668, 7680, 7692, 7710, 7716, 7722, 7728, 7734, 
7740, 7758, 7764, 7770, 
\\
7776, 7782, 7788, 7806, 7812, 7818, 7824, 7830, 7836, 7854, 7860, 7866, 7872, 
7884, 7902, 7914, 7920, 
\\
7926, 7932, 7950, 7956, 7968, 7980, 7998, 8004, 8010, 8016, 8022, 8028, 8046, 
8052, 8058, 8064, 8076, 
\\
8094, 8100, 8106, 8112, 8118, 8124, 8142, 8154, 8160, 8166, 8172, 8190, 8196, 
8202, 8208, 8214, 8220, 
\\
8238, 8244, 8250, 8256, 8262, 8268, 8286, 8298, 8304, 8316, 8334, 8340, 8346, 
8352, 8358, 8364, 8382, 
\\
8388, 8394, 8400, 8412, 8430, 8436, 8442, 8448, 8454, 8460, 8478, 8484, 8490, 
8496, 8502, 8508, 8526, 
\\
8532, 8538, 8544, 8550, 8556, 8574, 8580, 8586, 8592, 8598, 8604, 8622, 8628, 
8634, 8640, 8646, 8652, 
\\
8670, 8676, 8682, 8688, 8694, 8700, 8718, 8724, 8730, 8736, 8742, 8748, 8766, 
8772, 8778, 8784, 8790, 
\\
8796, 8814, 8820, 8826, 8832, 8844, 8862, 8868, 8874, 8880, 8886, 8892, 8910, 
8916, 8922, 8928, 8934, 
\\
8940, 8958, 8964, 8970, 8976, 8982, 8988, 9006, 9012, 9018, 9024, 9030, 9036, 
9054, 9060, 9066, 9072, 
\\
9078, 9084, 9102, 9108, 9114, 9120, 9126, 9132, 9150, 9156, 9162, 9168, 9174, 
9180, 9198, 9204, 9210, 
\\
9216, 9222, 9228, 9246, 9252, 9258, 9264, 9270, 9276, 9294, 9300, 9306, 9312, 
9318, 9324, 9342, 9348, 
\\
9354, 9360, 9366, 9372, 9390, 9396, 9402, 9408, 9414, 9420, 9438, 9444, 9450, 
9456, 9462, 9468, 9486, 
\\
9492, 9498, 9504, 9510, 9516, 9534, 9540, 9546, 9552, 9558, 9564, 9582, 9588, 
9594, 9600, 9606, 9612, 
\\
9630, 9636, 9642, 9648, 9654, 9660, 9678, 9684, 9690, 9696, 9702, 9708, 9726, 
9732, 9738, 9744, 9750, 
\\
9756, 9774, 9780, 9786, 9792, 9798, 9804, 9822, 9828, 9834, 9840, 9846, 9852, 
9870, 9876, 9882, 9888, 
\\
9894, 9900, 9918, 9924, 9930, 9936, 9942, 9948, 9966, 9972, 9978, 9984, 9990, 
9996, 10014, 10020, 
\\
10026, 10032, 10038, 10044, 10062, 10068, 10074, 10080, 
10086, 10092, 10110, 10116, 10122, 10128, 
\\
10134, 10140, 10158, 10164, 10170, 10176, 10182, 10188, 
10206, 10212, 10218, 10224, 10230, 10236, 
\\
10254, 10260, 10266, 10272, 10278, 10284, 10302, 10314, 
10320, 10326, 10332, 10350, 10356, 10362, 
\\
10368, 10374, 10380, 10398, 10404, 10410, 10416, 10422, 
10428, 10446, 10452, 10458, 10464, 10470, 
\\
10476, 10494, 10500, 10506, 10512, 10518, 10524, 10542, 
10548, 10560, 10566, 10572, 10590, 10596, 
\\
10602, 10608, 10614, 10620, 10638, 10644, 10650, 10656, 
10662, 10668, 10686, 10692, 10698, 10704, 
\\
10710, 10716, 10734, 10740, 10746, 10752, 10758, 10764, 
10782, 10788, 10794, 10800, 10806, 10812, 
\\
10830, 10836, 10842, 10848, 10860, 10878, 10884, 10890, 
10896, 10902, 10908, 10926, 10932, 10938, 
\\
10944, 10950, 10956, 10974, 10980, 10986, 10992, 10998, 
11004, 11022, 11034, 11040, 11046, 11052, 
\\
11070, 11076, 11082, 11088, 11094, 11100, 11118, 11124, 
11130, 11136, 11142, 11148, 11166, 11178, 
\\
11184, 11190, 11196, 11214, 11220, 11226, 11232, 11238, 
11244, 11262, 11268, 11280, 11286, 11292, 
\\
11310, 11322, 11328, 11334, 11340, 11358, 11364, 11370, 
11376, 11388, 11406, 11412, 11424, 11430, 
\\
11436, 11454, 11466, 11472, 11484, 11502, 11508, 11514, 
11520, 11532, 11550, 11562, 11568, 11580, 
\\
11598, 11610, 11616, 11622, 11646, 11652, 11658, 11664, 
11694, 11700, 11706, 11712, 11724, 11754, 
\\
11760, 11766, 11772, 11790, 11802, 11808, 11820, 11838, 
11844, 11856, 11868, 11886, 11892, 11898, 
\\
11904, 11934, 11940, 11946, 11952, 11964, 11982, 11994, 
12000, 12030, 12048, 12060, 12096, 12108, 
\\
12126, 12144, 12156, 12174, 12192, 12204, 12222, 12228, 
12240, 12270, 12288, 12336, 12342, 12348, 
\\
12366, 12378, 12414, 12432, 12444, 12462, 12516, 12588, 
12606, 12624, 12630, 12654, 12666, 12672, 
\\
12702, 12720, 12768, 12798, 12894, 12912, 12942, 12990, 
13020, 13086, 13098, 13200, 13278, 13308, 
\\
13536, 13566, 13596, 13662, 13854, 14142, 14250, 14460, 
14718 
\\
\noalign{\hrule height0.8pt}
\end{tabular}
}
\end{center}
\end{table}
%%%%%%%%%%%%%%%%%%%%%%%%%%%%%%%

We denote the $1532$ inequivalent bordered double circulant
doubly even self-dual $[96,48,16]$ codes by
$B_{96,i}$ $(i=1,2,\ldots,1532)$.
The first rows $r$ of $R'$ in~\eqref{bDCC} can be obtained from
\url{http://www.math.is.tohoku.ac.jp/~mharada/Paper/DCCb96.txt}.
By calculating the numbers of codewords of weight $16$,
we verified that
the $1532$ codes have $25$ distinct weight enumerators.
% The smallest and largest integers $a$ among the weight enumerators
% of $B_{96,i}$ $(i=1,2,\ldots,1532)$ are $6204$ and $13254$,
% respectively.
For each $a$,
the number $N_{a}$ of codes with the weight enumerator
is listed in Table~\ref{Tab:WD96}.
We verified by {\sc Magma}~\cite{Magma} that
the $1532$ codes have automorphism groups of order $94$.

%%%%%%%%%%%%%%%%%%%%%%%%%%%%%%%
\begin{table}[thbp]
\caption{Weight enumerators of bordered double circulant 
doubly even self-dual $[96,48,16]$ codes}
\label{Tab:WD96}
\begin{center}
%{\small
{\footnotesize
%{\scriptsize
\begin{tabular}{c|c|c|c|c|c|c|c|c|c}
\noalign{\hrule height0.8pt}
$a$ &  $N_a$ &$a$ &  $N_a$ &$a$ &  $N_a$ &$a$ &  $N_a$ &$a$ &  $N_a$ \\
\hline
$ 6204$ & $  1$&
$ 6768$ & $  2$&
$ 7050$ & $  5$&
$ 7332$ & $  8$&
$ 7614$ & $ 17$\\
$ 7896$ & $ 30$&
$ 8178$ & $ 72$&
$ 8460$ & $ 82$&
$ 8742$ & $116$&
$ 9024$ & $141$\\
$ 9306$ & $157$&
$ 9588$ & $197$&
$ 9870$ & $141$&
$10152$ & $160$&
$10434$ & $130$\\
$10716$ & $ 92$&
$10998$ & $ 74$&
$11280$ & $ 32$&
$11562$ & $ 31$&
$11844$ & $ 20$\\
$12126$ & $ 12$&
$12408$ & $  4$&
$12690$ & $  5$&
$12972$ & $  1$&
$13254$ & $  2$\\
\noalign{\hrule height0.8pt}
\end{tabular}
}
\end{center}
\end{table}
%%%%%%%%%%%%%%%%%%%%%%%%%%%%%%%

%%%%%%%%%%%%%%%%%%%%%
\section{Neighbors of double circulant self-dual $[96,48,16]$ codes}

Let $C$ be a double circulant self-dual $[96,48,16]$ code.
Let $M$ be a matrix whose rows are the codewords of weight $16$ in $C$.
We verified by {\sc Magma}~\cite{Magma} that
\[
(\rank M,\rank ( M \ \allone^T))=(47,48),
\]
for $C=C_{96,i}$ ($i=1,2,\ldots,49$) and
\[
(\rank M,\rank ( M \ \allone^T))=(48,49),
\]
for $C=P_{96,i}$ ($i=1,2,\ldots,4565$)
and $C=B_{96,i}$ ($i=1,2,\ldots,1532$).
% In addition, we verified by {\sc Magma} that
% the self-dual neighbors constructed by the above argument
% have minimum weight at most $16$.
By the method given in Section~\ref{sec:90nei},
we have the following results.

\begin{prop}
No double circulant self-dual $[96,48,16]$ code
has a self-dual $[96,48,d]$ neighbor with $d \ge 18$.
\end{prop}

It is still an open problem whether
an extremal doubly even self-dual code or
an extremal singly even self-dual code of length $96$ exists.

%%%%%%%%%%%%%%%%%%%%
\bigskip
\noindent
{\bf Acknowledgment.}
The authors would like to thank Alfred Wassermann
for his useful private communication.
This work was supported by JSPS KAKENHI Grant Number 15H03633.

%%%%%%%%%%%%%%%%%%%  References  %%%%%%%%%%%%%%%%%%%%%%%%

%%%%%%%%%%%%%%%%%%%%%%%%%%%%%%%%%


\begin{thebibliography}{99}

\bibitem{Magma}W. Bosma, J. Cannon and C. Playoust,
The Magma algebra system I: The user language,
{\sl J. Symbolic Comput.}
{\bf 24} (1997), 235--265.

\bibitem{CHK64}N. Chigira, M. Harada and M. Kitazume,
{Extremal self-dual codes of length $64$ through
neighbors and covering radii},
{\sl Des.\ Codes Cryptogr.}
{\bf 42} (2007), 93--101.

\bibitem{C-S}J.H. Conway and N.J.A. Sloane,
{A new upper bound on the minimal distance of self-dual codes},
{\sl IEEE Trans.\ Inform.\ Theory}
{\bf 36} (1990), 1319--1333.

\bibitem{DRZ}R. Dontcheva, R. Russeva and N. Ziapkov,
On the binary extremal codes of length 92,
{\sl Proc.\ Intern.\ Workshop Optimal Codes and Related Topics},
Sunny Beach, Bulgaria, 2001, pp.~53--58.

\bibitem{DGH}S.T.~Dougherty, T.A.~Gulliver and M.~Harada,
{Extremal binary self-dual codes,}
{\sl IEEE Trans.\ Inform.\ Theory}
{\bf 43} (1997), 2036--2047.

\bibitem{GO}P. Gaborit and A. Otmani,
Experimental constructions of self-dual codes,
{\sl Finite Fields Appl.}
{\bf  9}  (2003),  372--394.

\bibitem{GPSA}P. Gaborit, V. Pless, P. Sol\'e and O. Atkin,
Type II codes over $\FF_4$,
{\sl Finite Fields Appl.}
{\bf 8} (2002), 171--183.

\bibitem{GH-WE}T.A. Gulliver and M. Harada,
{Weight enumerators of double circulant codes and new extremal
self-dual codes},
{\sl Des.\ Codes and Cryptogr.}
{\bf 11} (1997), 141--150.

\bibitem{GH-DCC}T.A. Gulliver and M. Harada,
{Classification of extremal double circulant self-dual codes of
lengths $64$ to $72$},
{\sl Des.\ Codes and Cryptogr.}
{\bf 13} (1998), 257--269.

\bibitem{GH-DCC88}T.A. Gulliver and M. Harada,
{Classification of extremal double circulant self-dual codes of
lengths $74$--$88$},
{\sl Discrete Math.}
{\bf 306} (2006), 2064--2072.

\bibitem{GS}T.A. Gulliver and N. Senkevitch,
On a class of self-dual codes derived from quadratic residues,
{\sl IEEE Trans.\ Inform.\ Theory}
{\bf 45} (1999), 701--702.

\bibitem{HGK-DCC}M. Harada, T.A. Gulliver and  H. Kaneta,
{Classification of extremal double-circulant self-dual codes of
length up to $62$},
{\sl Discrete Math.}
{\bf 188} (1998), 127--136.

\bibitem{HKWY}M. Harada, M. Kiermaier, A. Wassermann and R. Yorgova,
New binary singly even self-dual codes,
{\sl IEEE Trans.\ Inform.\ Theory}
{\bf  56}  (2010),  1612--1617.

\bibitem{HN}M. Harada and T. Nishimura,
An extremal singly even self-dual code of length 88,
{\sl Adv.\ Math.\ Commun.}
{\bf 1} (2007), 261--267.

%\bibitem{Karlin}M. Karlin,
%New binary coding results by circulants,
%{\sl IEEE Trans.\ Inform.\ Theory},
%{\bf 15} (1969), 81--92.

\bibitem{Mallows-Sloane}C.L.~Mallows and N.J.A.~Sloane,
{An upper bound for self-dual codes},
{\sl Inform.\ Control}
{\bf 22} (1973), 188--200.

\bibitem{Moore}E.H. Moore,
{\sl Double circulant codes and related structures},
Ph.D. dissertation, Dartmouth College, Hanover, NH, (1976).

% \bibitem{TTHAA}C. Tjhai, M. Tomlinson, R. Horan, M. Ahmed and M. Ambroze,
% Some results on the weight distributions of the binary
% double-circulant codes based on primes,
% {\sl Proc.\ IEEE Singapore Inter.\ Conference Commun. \ Systems},
% Singapore, 2006.
% http://ieeexplore.ieee.org/document/4085726/
% I did not know the page numbers.

\bibitem{Y92}R. Yorgova,
Binary self-dual extremal codes of length 92,
{\sl Proc.\ IEEE Intern.\ Symp.\ Inform.\ Theory},
Seattle, WA, 2006, pp.~1292--1295.

\bibitem{YW}R. Yorgova and A. Wassermann,
Binary self-dual codes with automorphisms of order 23,
{\sl Des.\ Codes Cryptogr.}
{\bf 48} (2008), 155--164.

\end{thebibliography}
\end{document}